\documentclass[3p]{elsarticle}

\makeatletter
\def\ps@pprintTitle{
	\let\@oddhead\@empty
	\let\@evenhead\@empty
	\let\@oddfoot\@empty
	\let\@evenfoot\@oddfoot
}
\makeatother

\usepackage{amssymb,amsmath,mathtools}

\usepackage{graphicx}
\usepackage{cancel}
\usepackage[ruled,vlined]{algorithm2e}
\usepackage{hyperref}
\usepackage[colorinlistoftodos,textwidth=4cm,shadow]{todonotes}

\newcounter{Igor}

\newcounter{Vasileios}

\date{}
\title {\textbf{Rigid Clumps in the \textit{\textit{MercuryDPM}} Particle Dynamics Code}}

\author[1]{Igor Ostanin\footnote{Corresponding author, e-mail:i.ostanin@utwente.nl.}}
\author[2,3]{Vasileios Angelidakis}
\author[1]{Timo Plath}
\author[1]{Sahar Pourandi}
\author[1]{Anthony Thornton} 
\author[1]{Thomas Weinhart}

\address[1]{Multi-Scale Mechanics (MSM), Faculty of Engineering Technology, MESA+, University of Twente, P.O. Box 217, 7500 AE Enschede, The Netherlands.}
\address[2]{Institute for Multiscale Simulation, Friedrich-Alexander-Universit\"at Erlangen-N\"urnberg, 91058, Erlangen, Germany }
\address[3]{ School of Engineering, Newcastle University, NE1 7RU, Newcastle upon Tyne, United Kingdom }

\begin{document}

\begin{abstract}
Discrete particle simulations have become the standard in science and industrial applications exploring the properties of particulate systems. Most of such simulations rely on the concept of interacting spherical particles to describe the properties of particulates, although, the correct representation of the nonspherical particle shape is crucial for a number of applications. In this work we describe the implementation of clumps, i.e. assemblies of rigidly connected spherical particles, which can approximate given nonspherical shapes, within the \textit{MercuryDPM} particle dynamics code. \textit{MercuryDPM} contact detection algorithm is particularly efficient for polydisperse particle systems, which is essential for multilevel clumps approximating complex surfaces. We employ the existing open-source \texttt{CLUMP} library to generate clump particles. We detail the pre-processing tools providing necessary initial data, as well as the necessary adjustments of the algorithms of contact detection, collision/migration and numerical time integration. The capabilities of our implementation are illustrated for a variety of examples.
	
\end{abstract}

\maketitle



\section{Introduction}

\subsection{Overview and scope}

Rigid assemblies of spherical particles \cite{Gallas1993, Nolan1995} are an important tool to simulate materials consisting of particles of irregular shapes with the discrete element method (DEM). The alternative approaches, that are often employed to model non-spherical particles in DEM \cite{Govender2016, Podlozhnyuk2017},  have certain limitations -- polyhedral particle shapes \cite{Govender2016} lead to difficulties in generalization of a wide set of well-established contact models for spherical particles, while superquadrics \cite{Podlozhnyuk2017} do not offer sufficiently general particle shape representation toolkit. As a result, almost all modern commercial DEM codes, e.g. EDEM \cite{EDEM} or PFC \cite{PFC}, include functionality to model rigid assemblies of spherical particles. 

However, as will be demonstrated below, the implementation of rigid clumps in DEM introduces ambiguities that are hard to interpret when the source code and exact implementation details are unavailable. We seek to fill this gap, presenting fully functional, well-documented and completely open source implementation of rigid particle assemblies within the \textit{MercuryDPM} \cite{Mercury2020} particle dynamics code, utilizing CLUMP library \cite{clump2021} for particle generation.

Below we provide a brief overview of the \textit{MercuryDPM} particle dynamics engine and discuss the notion of a \textit{rigid clump} -- a rigid assembly of spherical particles -- as it will be used in this paper. In the following sections we will take a closer look at the necessary theoretical background, the implementation details and the examples of using rigid clumps in numerical simulations with \textit{MercuryDPM}.      

\subsection{\textit{MercuryDPM} particle dynamics code}

\textit{MercuryDPM} \cite{MDPMBitbucket} is an open-source realization of DEM. It is mainly used to simulate granular particles -- collections of discrete particles that can be found in many natural and artificial settings. Examples include snow, sand, soil, coffee, rice, coal, pharmaceutical tablets, catalysts, and animal feed. Understanding the behavior of such materials is crucial for industries like pharmaceuticals, mining, food processing, and manufacturing.  

The development of the code started in 2009 at the University of Twente, and since then it has grown into a large framework with a wide open-source community of academic and industrial users. The core development team is still located at the University of Twente. \textit{MercuryDPM} is a versatile, object-oriented C++ code that is built and tested using the capabilities of \texttt{cmake/ctest}. 

The code possesses three primary features enabling it to simulate complex industrial and natural scenarios: (i) the flexible implementation allowing complex walls and boundary conditions; (ii) the analysis toolkit, able to extract the most relevant information from the large amount of data generated by these simulations, (iii) the advanced contact detection scheme that makes \textit{MercuryDPM} particularly efficient for highly polydisperse particle systems; \cite{Ogarko2012, Mercury2020}. The latter feature is particularly interesting in a context of simulating clumps, since fine representation of shape of a non-spherical particle often requires highly polydisperse clumps. 

\textit{MercuryDPM} normally operates with spherical particles (discrete elements), characterized by the mass, radius, position, velocity and angular velocity. Also \textit{MercuryDPM} offers support of superquadric particles \cite{Mercury2020}. The Velocity Verlet time integration algorithm is utilised to update the positions of each particle, while the forward Euler algorithm is employed for particle rotations. Particle interactions are governed by wide variety of contact models which describe physical laws to compute the normal and tangential forces resulting from particle's contacts.

\subsection{Rigid clumps}

By \textit{rigid clump} (or just \textit{clump}) we will imply an aggregate of $N$ rigid spherical particles of a given density, that are rigidly linked to each other at a given relative translational and rotational positions (Fig. 1). The constituent particles of a clump will be referred to as \textit{pebbles}. The clump is a \textit{rigid body} possessing $6$ degrees of freedom. Therefore, in 3D, the number of constraints that are implicitly introduced on relative translational and rotational positions of particles is $6 (N-1)$.

The pebbles may (or may not) have overlaps, introducing volumes within a clump that belong to more than one pebble. It is therefore impossible to algebraically sum up the inertia of the clump over pebbles for a system of overlapping pebbles representing a complex-shaped particle. Our approaches to computing the inertia of clumps are discussed below. 

\begin{figure}
	\begin{center}
		\includegraphics[width=10.0cm]{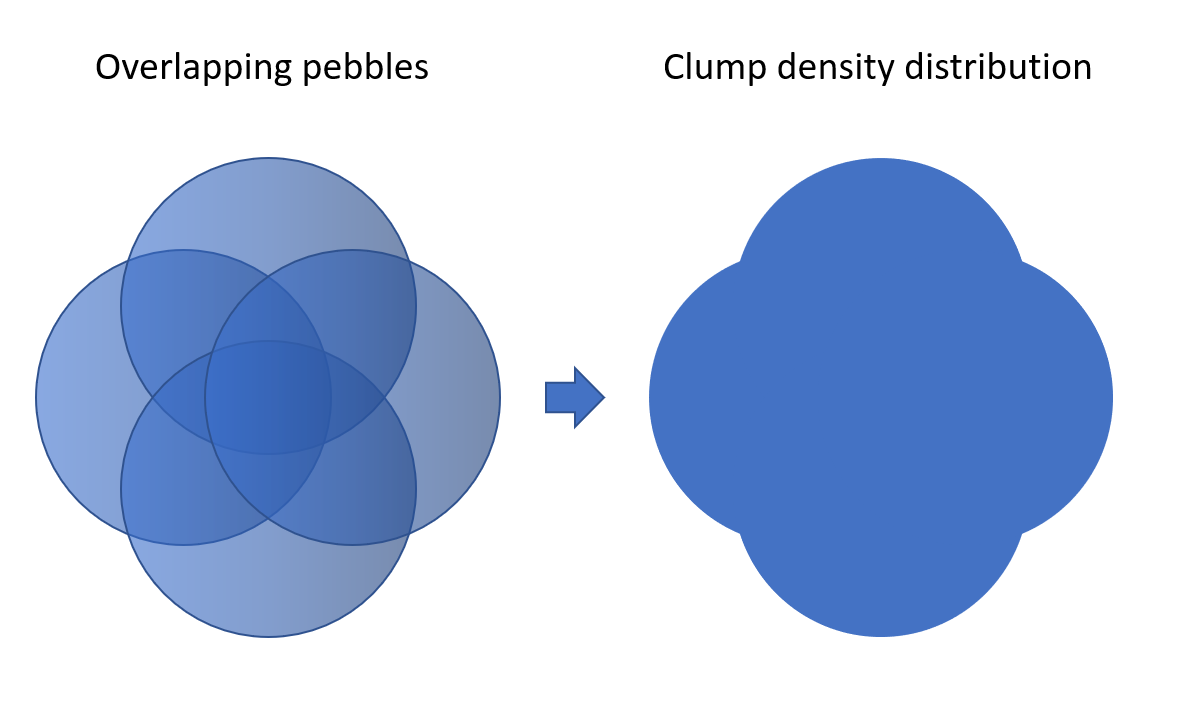}
		\protect\caption{Rigid clump and its inertial properties -- conceptual illustration}
	\end{center}
\end{figure}

Our implementation builds on the multispheres featured in the earlier versions of the code (see section 6.2 in \cite{Mercury2020}). However, the functional and performance of the implementation have been significantly expanded and improved via incorporation of multiple new features, architecture improvements and bugfixes. The new implementation allows to address a wide class of problems that previously remained unavailable - large simulation model sizes, arbitrarily complex clump geometries, complex (e.g. moving periodic) boundary conditions etc.


\section{Clump geometry generation with the \texttt{CLUMP} software}

\texttt{CLUMP} software \cite{clump2021} has been developed recently to address the problem of automatic generation of rigid clump particles by approximation of polyhedral shapes. MercuryDPM provides the necessary interface to use \texttt{CLUMP}-generated particles in DEM simulations. This section offers an overview of the main features of \texttt{CLUMP} and underlying clump generation methodologies.

The open-source \texttt{CLUMP} software (Code Library to generate Universal Multi-sphere Particles) \cite{clump2021} is used to create clump representations of irregular particle geometries. This software takes as input shape/imaging data of various types, such as point clouds, surface meshes (e.g. in the form of stereolithography/stl files), tetrahedral meshes and labelled three dimensional images derived via Computed Tomography. Utilising this input, three clump generation methods are implemented in the software to create clump representations of them, proposed in \cite{favier1999shape}, \cite{ferellec2010method} and \cite{clump2021}.

The method of \cite{favier1999shape}, one of the historically first clump generation methods, is implemented in the software to generate clumps of axisymmetric bodies. Although the original paper introducing the method \cite{favier1999shape} does not delineate a way of generating clumps of real particles, the implementation in \texttt{CLUMP} offers the capability of achieving this via the following steps: a particle geometry is loaded from imaging data and its inertia tensor is calculated; the principal inertia values and principal directions (PDs) are determined and the particle is oriented to its PDs; then a user-defined number of spheres are generated along the longest particle dimension, the size of which is decided so as to approximate the shape of the input particle; last, the clump is oriented back to the original orientation of the input particle. This method can generate elongated and compact particles of limited elongation, but cannot generate particles with pronounced flat features.

For irregular particles that do not display axisymmetric features, the method described in \cite{ferellec2010method} is an efficient method to generate clumps based on a triangulated mesh representation of the particle surface (i.e. made of vertices and triangular faces). The method first calculates the normal vectors of each vertex as the average of the adjacent face normal vectors; then, a random vertex is selected and a tangent sphere is grown internally within the particle, until it intersects one of the other particle vertices; the process is repeated until a sought number of spheres is generated. If imaging data are given in a different format, e.g. via Computed Tomography, this is handled internally within \texttt{CLUMP}, via transformation of the data to a surface mesh. The simplicity of the method makes it appealing and computationally efficient, but the random selection of vertices can lead to inadequate clumps for small numbers of spheres per particle. In such cases, for the same number of spheres the algorithm generates clumps of vastly different characteristics, as there is no rationale behind the random selection of vertices. As a result, for these cases there is no correlation between employed number of spheres and achieved morphological fidelity. However, if a large amount of spheres is considered computationally affordable by the modeller (e.g. in \cite{ferellec2010method} up to $5500$ spheres were considered), this method generates clumps with reduced artificial surface roughness, as reported in \cite{ferellec2010method}.

A new clump generation technique was recently proposed as part of \texttt{CLUMP} \cite{clump2021}, which relies on the Euclidean transform of three-dimensional images. A particle shape is either imported directly from binarized (or labelled) images, or transformed into a three-dimensional image from other data types (e.g. from surface mesh data); the Euclidean transform of the image is calculated, and the maximum value of the transform determines the location and radius of the largest possible inscribed sphere that fits in the particle. This sphere is considered as the first sphere, the voxels corresponding to a percentage of this sphere are deactivated from the original image, leading to a residual image (original minus a percentage of the sphere voxels); then, the Euclidean transform of this new residual image is used to calculate the next sphere; the process is repeated until a user-defined required number of spheres is generated or if a user-defined minimum radius is achieved. This technique has the clear advantage that each new sphere is generated at the position where the mass of the particle is least represented, thus creating a clear correlation between the number of spheres (a descriptor associated to computational cost) and the achieved morphological similarity (a descriptor of simulation fidelity). With this method, each sphere is of equal or smaller size to its previous one, and so particle generation is performed in a systematic and predictable way. If all the voxels of a sphere are deactivated after each iteration of the method, the method results in clusters of non-overlapping spheres, while if only a percentage of each sphere is deactivated, clumps of overlapping spheres are generated, as delineated in \cite{clump2021}. The drawback of the method is its high cost in terms of memory consumption (though still manageable even for a regular desktop computer).

Choosing the optimal or preferred particle generation technique lies with the user, as different applications and different particle types pose different requirements in terms of the employed particle characteristics. In terms of efficiency, all of the aforementioned techniques perform well, mainly due to their algorithmic simplicity, allowing for the generation of several hundred particles within few minutes, for input imaging data of reasonable resolution and size.


\section{Rigid clumps in \textit{MercuryDPM}}

\subsection{General organization}

The rigid clump functional in \textit{MercuryDPM} is currently implemented as a multilevel structure. The logic of unification of pebbles in the clump, as well as the algorithms of time integration are implemented in the class \texttt{./Kernel/particles/ClumpParticle.h/cc} inherited from an abstract nonspherical particle class \texttt{./Kernel/particles/NonSphericalParticle.h/cc}, that, in turn is inherited from the base particle class \texttt{./Kernel/particles/BaseParticle.h/cc}. It is expected that the functions inherent to all types of nonspherical particles (e.g. rigid dynamics time integration) in the future will be located in the class \texttt{./Kernel/particles/NonSphericalParticle.h/cc}.
   
The \texttt{CLUMP} software, described above, is used to generate positions and radii of pebbles that describe the given nonspherical shape. The \texttt{CLUMP} tool provides pebble data, which, along with the optionally provided initial \texttt{stl} format shape of the clump, constitute an input of \texttt{MClump} pre-processing tool (part of \textit{MercuryDPM}, cite \cite{MClump}). Alternatively, the pebble data for \texttt{MClump} can be generated manually.  

\texttt{MClump} centers and rotates the clump, aligning its
principal axes with the global Cartesian axes, and computes clump’s inertia using the prescribed algorithm (summation over pebbles, summation over voxels, summation over tetrahedrons using \texttt{stl} representation - see the description below. Fig. 2. details modes of work of \texttt{MClump}. In the first mode, \texttt{MClump} imports list of pebbles and then does all the necessary computations (center of mass (COM), volume, tensor of inertia (TOI), principal directions) based on summation over pebbles, as discussed in Subsection 3.3.1. In the second mode, \texttt{MClump} imports list of pebbles, but performs inertia computations on the voxel grid, excluding extra contributions of pebble's overlaps (Subsection 3.3.2). In the third mode, \texttt{MClump} receives the triangulated surface of a nonspherical particle, as well as its clumped sphere approximation generated by the external tool (\texttt{CLUMP} library), and computes the necessary properties (Subsection 3.3.3).   

Headers for the driver files

\texttt{./Drivers/Clump/ClumpHeaders/ClumpIO.h},

\texttt{./Drivers/Clump/ClumpHeaders/Mercury3DClump.h}, 

introduce necessary features and modifications of the \textit{MercuryDPM} virtual members, enabling clump dynamics, namely:

\begin{itemize}
    \item The modifications of the \textit{MercuryDPM} engine, changing the logic of application of contact forces and moments, as well as the external forces (e.g. gravity).

    \item The adjustment of the logic of interaction of the clump and its pebbles with the periodic boundary.
    
    \item The import tool, that loads the all data of available clump instances, including clump volume, TOI and the list of pebbles. 
    
    \item Clump distribution generation functions, that create distributions of non-overlapping rotated clumps in a given spatial domain.
\end{itemize}

Driver files (compiled simulation descriptions, see \cite{Mercury2020} for details) utilize these tools to load the list of clump instances generated by \texttt{MClump}, and, using them, generate necessary distributions of clumps and compute their dynamics.

\begin{figure}
	\begin{center}
		\includegraphics[width=14.5cm]{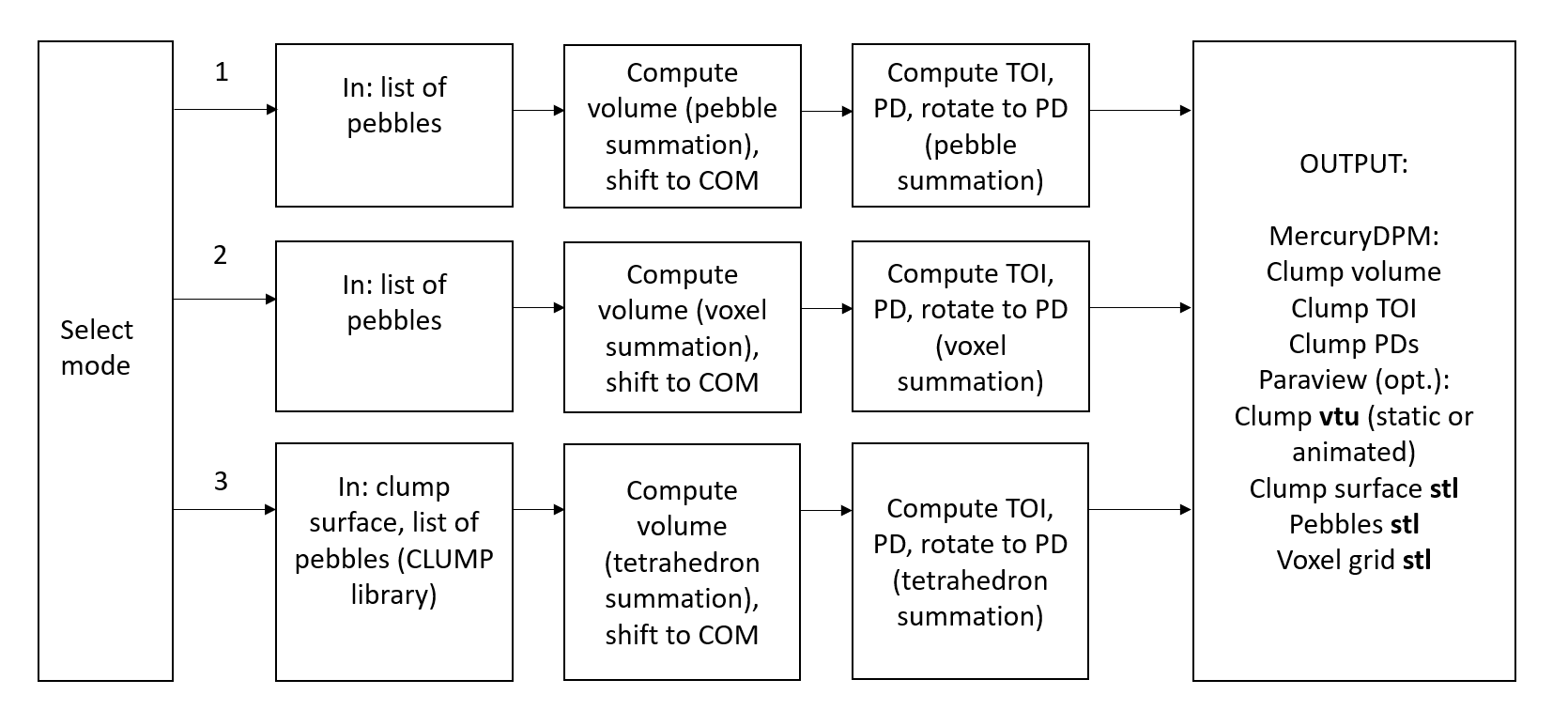}
		\protect\caption{Modes of operation of MClump tool.}
	\end{center}
\end{figure}

\subsection{Clump creation logic}

The unification of particles into rigid clumps occurs by assigning to every particle instance the role of either a ``clump'' particle or a ``pebble'' particle. Specifically, every instance of BaseParticle class has Boolean attributes (flags) \texttt{isClump} and \texttt{isPebble}. The ``pebble'' instances have \texttt{isClump = False}, \texttt{isPebble = True}. All the ``clump'' (container) instances have \texttt{isClump = True}, \texttt{isPebble = False}. Regular spherical particles have \texttt{isClump = False}, \texttt{isPebble = False}. Depending on the flags, these three \footnote{The ${4-}^{th}$ combination (\texttt{isClump = True}, \texttt{isPebble = True}) is explicitly excluded in the relevant particle class methods.} types of particles have different behavior in contact detection, migration over boundaries etc. Namely, for the clump particle the interactions are treated at the pebble level, while time integration of motion occurs at the clump level. Stiff pebble-pebble interactions are assumed, so that clump-clump contact is always represented with a single pebble-pebble contact. In case of (non-physical) multiple contacts between the close pebbles of interacting clumps, the corresponding effective increase in stiffness should be treated \cite{Kodam2009}. The motion of pebbles is prescribed according to translation and rotation of the corresponding clump. Clumps and pebbles have some other differences in behavior, e.g. in a context of interaction with periodic boundaries -- see the discussion below. 


\subsection{Computing inertial properties of a clump}

Defining inertial properties of a clump is a non-trivial problem. The analytical treatment is possible in case of absent overlaps (direct summation over pebbles, as implemented earlier \cite{Mercury2020}), and overlaps between no more than two spherical pebbles (summation over pebbles and subtraction of "cap" segments, \cite{Parteli_2013}). In our implementation, we use three different approaches to compute mass and TOI of complex shape particles: summation over the pebbles, summation over the voxels and summation over the tetrahedrons. Fig. 3. gives the qualitative idea about these representations of the volume of a non-spherical particle. Let us take a closer look at each of these approaches.

\subsubsection{Summation over pebbles}

This method of computation works if the pebbles do not overlap or we presume that the inertial properties of a clump are defined by the total mass of the pebbles. In this case the total mass and TOI can be directly summed over the spherical pebbles using mass conservation and Steiner's theorem. Given the density of pebbles $\rho$, their radii $r_j$ and positions in Cartesian system $\mathbf{x}_j = (x_j, y_j, z_j)$, we first find the mass of the clump and the position of the center of mass:

\begin{equation} \label{eq1}
    M = \sum m_j = \sum \frac{4}{3} \pi r_j^3 \rho
\end{equation}

\begin{equation} \label{eq2}
    \mathbf{x}_c = \frac{1}{M} \sum m_j \mathbf{x}_j
\end{equation}

At the next step, we shift the center of the coordinate system to the center of mass:

\begin{equation} \label{eq3}
     \mathbf{x}_j := \mathbf{x}_j - \mathbf{x}_c
\end{equation}

then we compute the TOI by summing over pebbles:

\begin{equation} \label{eq4}
    \mathbf{I} = \sum \mathbf{I}_j
\end{equation}

\begin{equation} \label{eq5}
	\mathbf{I}_j = m_j \begin{pmatrix}
		\frac{2}{5} r_j^2 + y_j^2 + z_j^2 & -x_j y_j & -x_j z_j\\
		-x_j y_j & \frac{2}{5} r_j^2 + x_j^2 + z_j^2 & -y_j z_j\\
		-x_j z_j & -y_j z_j & \frac{2}{5} r_j^2 + x_j^2 + y_j^2\\
	\end{pmatrix}
\end{equation}

Given the above-mentioned assumptions, this method is precise.

\begin{figure}
	\begin{center}
		\includegraphics[width=14.5cm]{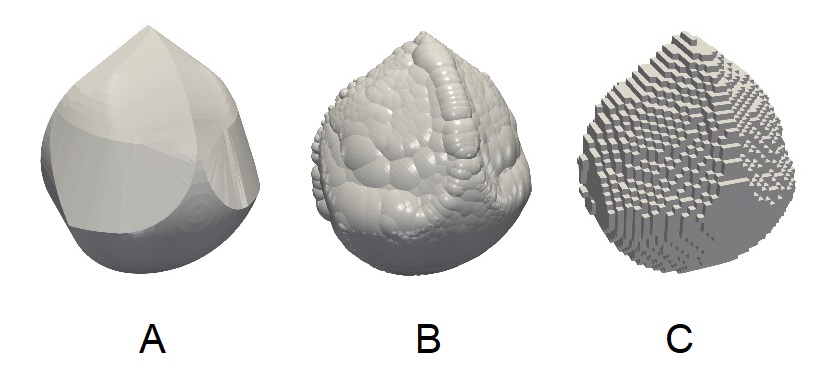}
		\protect\caption{Representation of a non-spherical shape as (A) triangulated surface, (B) rigid clump of spherical particles, (C) 3D array of voxels.}
	\end{center}
\end{figure}

\subsubsection{Summation over voxels}

In case if pebbles overlap and clump was not generated from a triangulated surface, we use voxel discretization to compute mass and TOI of a clump. The bounding box encapsulating every point of the clump is expanded to the cubic box $(x_b, x_b + Nd, y_b, y_b + Nd, z_b, z_b + Nd)$ of a minimal size, which is split into cubic voxels of side $d$, defined by the specified number of voxels $N$ along the side of a bounding box. Then the mask $\mathcal{M}(m,n,k)$ is introduced: $\mathcal{M}(m,n,k) = 1$ if the center of the voxel $m,n,k$ is inside of at least one pebble, and $\mathcal{M}(m,n,k) = 0$ otherwise. The coordinates of the center of the voxel are found as $\mathbf{x}(m, n, k) = (x_b + d(m+0.5), y_b + d(n+0.5), z_b + d(k+0.5)) $. Then the mass and the COM is computed as

\begin{equation} \label{eq6}
    M = \sum_m \sum_n \sum_k \mathcal{M}(m,n,k)  \rho d^3   
\end{equation}

\begin{equation} \label{eq7}
    \mathbf{x}_c = \frac{1}{M} \sum_m \sum_n \sum_k \mathcal{M}(m,n,k) \mathbf{x}(m,n,k) \rho d^3   
\end{equation}

Next, we shift the center of the coordinate system to the center of mass:

\begin{equation} \label{eq8}
     \mathbf{x}_j := \mathbf{x}_j - \mathbf{x}_c
\end{equation}

then we compute the TOI by summing over voxels:

\begin{equation} \label{eq9}
    \mathbf{I} = \sum_m \sum_n \sum_k \mathbf{I}_{mnk}
\end{equation}

\begin{equation} \label{eq10}
	\mathbf{I}_{mnk} = \rho d^3 \begin{pmatrix}
		y_{mnk}^2 + z_{mnk}^2 & -x_{mnk} y_{mnk} & -x_{mnk} z_{mnk}\\
		-x_{mnk} y_{mnk} & x_{mnk}^2 + z_{mnk}^2 & -y_{mnk} z_{mnk}\\
		-x_{mnk} z_{mnk} & -y_{mnk} z_{mnk} &   x_{mnk}^2 + y_{mnk}^2\\
	\end{pmatrix}
\end{equation}

The precision of such estimate depends on the chosen resolution and the complexity of the shape. The method requires brute-force summation over ~$10^6-10^9$ voxels, therefore, pre-computation might take some time.

\subsubsection{Summation over tetrahedrons}

If the clump is generated by approximation of known triangulated surface, one can use the latter for an explicit calculation of the TOI \cite{TOI_tetra}. In this case the TOI is computed by the analytical summation over tetrahedrons.

The COM of a tetrahedron $j$ with the vertices $[\mathbf{a}^1_j, \mathbf{a}^2_j, \mathbf{a}^3_j, \mathbf{a}^4_j]$ is given by:

\begin{equation} \label{eq11}	
\mathbf{c}_j = \frac{\mathbf{a}^1_j + \mathbf{a}^2_j + \mathbf{a}^3_j + \mathbf{a}^4_j}{4}
\end{equation}

Volume of a tetrahedron $j$ is given by

\begin{equation} \label{eq12}	
	V_j = \frac{1}{6} \begin{vmatrix}
		(a^1_j)_x & (a^1_j)_y & (a^1_j)_z & 1\\
		(a^2_j)_x & (a^2_j)_y & (a^2_j)_z & 1\\
		(a^3_j)_x & (a^3_j)_y & (a^3_j)_z & 1\\
		(a^4_j)_x & (a^4_j)_y & (a^4_j)_z & 1\\ 
	\end{vmatrix}
\end{equation}.

Here the volume $V_j$ comes with the the sign that depends on whether the normal $(\mathbf{a}_3 - \mathbf{a}_2) \times (\mathbf{a}_4 - \mathbf{a}_2)$ is directed into the half-space containing the vertex $\mathbf{a}_1$ (negative sign) or vise versa.

Given arbitrary volume bounded by a triangulated surface $\Gamma$ consisting of a set of triangles $\mathbf{s}_j$, and an arbitrary point $\mathbf{O}$, one can compute the COM of the volume as:

\begin{equation} \label{eq13}	
	\mathbf{x}_c = \sum \mathbf{c}_j V_j
\end{equation}
   
where $V_j$ and $\mathbf{c}_j$ are the volume and COM of the tetrahedron $[\mathbf{a}_1, \mathbf{a}_2, \mathbf{a}_3, \mathbf{a}_4] = [\mathbf{O}, \mathbf{s}^1_j, \mathbf{s}^2_j, \mathbf{s}^3_j]$

Similarly to alternative approaches, we shift the coordinate system to match its origin with the computed clump's COM:

\begin{equation} \label{eq14}
     \mathbf{x}_j := \mathbf{x}_j - \mathbf{x}_c
\end{equation}

One can further compute mass and TOI of the body with respect to its COM as the sum of masses and moments of inertia of constituent tetrahedrons:

\begin{equation} \label{eq15}
	\begin{split}
		M = \rho \sum V_j \\ 	
		\mathbf{I} = \sum \mathbf{I}_j
	\end{split}
\end{equation}

Here the TOI of a tetrahedron with respect to its first vertex $\mathbf{a}_1$, corresponding to the origin of the coordinate system and clump's COM, is computed according to \cite{TOI_tetra}:

\begin{equation} \label{eq16}	
	\mathbf{I}_j = \rho \begin{pmatrix}
		a & -c' & -b'\\
		-c' & b & -a'\\
		-b' & -a' & c\\
		\end{pmatrix}
\end{equation}

where

\begin{align} \label{eq17}	
	a &= \int_{D} (y^2+ z^2) dD, &
	b &= \int_{D} (x^2+ z^2) dD, &
	c &= \int_{D} (x^2+ y^2) dD, \\ \nonumber
	a' &= \int_{D} yz dD, &
	b' &= \int_{D} xz dD, &
	c' &= \int_{D} xy dD, \\ \nonumber
\end{align}

where $D$ is the tetrahedral domain.

It worth noting that the paper \cite{TOI_tetra} has a known \cite{Vasileios_bugfix} error, that is fixed in (\ref{eq16}): components $b'$ and $c'$ are erroneously swapped there.  

Denoting $[\mathbf{a}_1, \mathbf{a}_2, \mathbf{a}_3, \mathbf{a}_4] =$ $[(x_1, y_1, z_1)$, $(x_2, y_2, z_2)$, $(x_3, y_3, z_3)$, $(x_4, y_4, z_4)]$, the integrals \ref{eq8} are solved explicitly as:

\begin{equation} \label{eq18}
\begin{split} 
	a &=V_j( y_1^2 + y_1y_2 + y_2^2 + y_1y_3 + y_2y_3 +
	y_3^2 + y_1y_4 + y_2y_4 + y_3y_4 + y_4^2 \\ & + z_1^2 + z_1z_2 + z_2^2 + z_1z_3 + z_2z_3 + z_3^2 + z_1z_4 + z_2z_4 + z_3z_4 + z_4^2)/10	
\end{split}
\end{equation}

\begin{equation} \label{eq19}
	\begin{split} 
	b &= V_j(x_1^2 + x_1x_2 + x_2^2 + x_1x_3 + x_2x_3 +x_3^2 + x_1x_4 + x_2x_4 + x_3x_4 + x_4^2 \\ & + z_1^2 + z_1z_2 + z_2^2 + z_1z_3 + z_2z_3 + z_3^2 + z_1z_4 + z_2z_4 + z_3z_4 + z_4^2)/10 	
\end{split}
\end{equation}

\begin{equation} \label{eq20}
	\begin{split} 
	c &= V_j(x_1^2 + x_1x_2 + x_2^2 + x_1x_3 + x_2x_3 + x_3^2 + x_1x_4 + x_2x_4 + x_3x_4 +  x_4^2 \\ & + y_1^2 + y_1y_2 + y_2^2 + y_1y_3 + y_2y_3 + y_3^2 + y_1y_4 + y_2y_4 + y_3y_4 + y_4^2)/10	
\end{split}
\end{equation}

\begin{equation} \label{eq21}
	\begin{split} 
	a' &= V_j(2y_1z_1 + y_2z_1 + y_3z_1 + y_4z_1 + y_1z_2 + 2y_2z_2 + y_3z_2 + y_4z_2 + y_1z_3 \\ & + y_2z_3 + 2y_3z_3 + y_4z_3 + y_1z_4 + y_2z_4 + y_3z_4 + 2y_4z_4)/20 	
\end{split}
\end{equation}

\begin{equation} \label{eq22}
	\begin{split} 
	b' &= V_j(2x_1z_1 + x_2z_1 + x_3z_1 + x_4z_1 + x_1z_2 + 2x_2z_2 + x_3z_2 + x_4z_2 + x_1z_3 \\& + x_2z_3 + 2x_3z_3 + x_4z_3 + x_1z_4 + x_2z_4 + x_3z_4 + 2x_4z_4)/20 	
\end{split}
\end{equation}

\begin{equation} \label{eq23}
	\begin{split} 
	c' &= V_j(2x_1y_1 + x_2y_1 + x_3y_1 + x_4y_1 + x_1y_2 + 2x_2y_2 + x_3y_2 + x_4y_2 + x_1y_3 \\& + x_2y_3 + 2x_3y_3 + x_4y_3 + x_1y_4 + x_2y_4 + x_3y_4 + 2x_4y_4)/20	
\end{split}
\end{equation}

This method gives precise TOI of the initial triangulated surface. It worth noting that the formulae (\ref{eq13}, \ref{eq15}) work for rather complex (non-convex, multiply connected) domains: if the absolute volume of tetrahedrons is higher than the volume of a body, the extra volume is swept twice with tetrahedrons of positive and negative volume computed according to \ref{eq12}, which results in correct values for body's total volume, mass, COM and TOI. The examples section gives the comparison of the methods to compute inertial properties of a clump used in our work.

\subsection{Computing the clump's PDs}
\medskip

Principal axes of inertia $\mathbf{e}_1, \mathbf{e}_2, \mathbf{e}_3$ are found as eigenvectors of $\mathbf{I}$:

\begin{equation} \label{eq24}
	\mathbf{I} \mathbf{e_i} = \lambda_i \mathbf{e_i}
\end{equation}

PDs are assured to form the \textit{right-handed} Cartesian basis.

Once the PDs of the clump's TOI are computed, the clump instance is rotated to align its PDs with the Cartesian axes:

\begin{equation} \label{eq25}
	\mathbf{x} := \mathbf{Q} \mathbf{x} 
\end{equation}

\begin{equation} \label{eq26}
	\mathbf{I} := \mathbf{Q}^T \mathbf{I} \mathbf{Q}
\end{equation}

where $\mathbf{Q}$ is the rotation matrix defined as

\begin{equation} \label{eq27}	
	\mathbf{Q} = \begin{pmatrix}
		\mathbf{n}_1\mathbf{e}_1 & \mathbf{n}_2\mathbf{e}_1 & \mathbf{n}_3\mathbf{e}_1\\
		\mathbf{n}_1\mathbf{e}_2 & \mathbf{n}_2\mathbf{e}_2 & \mathbf{n}_3\mathbf{e}_2\\
		\mathbf{n}_1\mathbf{e}_3 & \mathbf{n}_2\mathbf{e}_3 & \mathbf{n}_3\mathbf{e}_3\\
		\end{pmatrix}
\end{equation}

where $\mathbf{n}_i$ are the orths of global Cartesian coordinate system, and $\mathbf{e}_i$ are orths of clump's eigendirections.

\subsection{Equations of motion of a rigid clump}

Once we have procedures that compute overall force $\mathbf{F}$ and moment $\mathbf{M}$ acting on the clump, we can solve the equations of motion using one of the schemes of numerical integration. For translational motion of a clump, we use the velocity Verlet algorithm that does not differ from the one employed for spherical particles, given that the particle mass is the mass of a clump. Below we consider the equations of motion for rotational degrees of freedom.

In the case when the TOI is non-spherical (the principal moments of inertia are not equal) the rotational dynamics is described by Euler equations:

\begin{equation} \label{eq28}	
	I_{ii} \dot \omega_i - I_{ij} \dot \omega_j + \epsilon_{ijk} \omega_j (I_{kk} \omega_k - I_{kl} \omega_l)) = M_i; (i \neq j, l \neq k)  	
\end{equation}

The non-spherical TOI $I_{ij}$ is computed based on one of the algorithms discussed above.

\subsection{Time integration of the EoM of a rigid clump}
\medskip

The time integration scheme used in our code utilizes a leap-frog algorithm of the time integration of the notion of non-spherical particle, similar to one utilized in PFC 4.0 \cite{pfc2008}. We track the orientation in the shape of rotation matrix $Q$ that in used to reconstruct the current orientation of local coordinate system and the positions of pebbles.  
The equation (\ref{eq27}) is solved using finite difference procedure of the second order, computing angular velocities $\omega_j$ at mid-intervals $t + \Delta t/2$, and all other quantities at primary intervals $t + \Delta t$. The equation (\ref{eq27}) can be re-written in the matrix form as

\begin{equation} \label{eq29}
	\begin{split}
		\mathbf{M} - \mathbf{W} &= \mathbf{I}  \dot \omega\\
		M &= \begin{pmatrix}
		M_1\\
		M_2\\
		M_3\\
		\end{pmatrix} \\
		W &= \begin{pmatrix}
		(I_{33} - I_{22})\omega_2 \omega_3 
		+ I_{23} \omega_3 \omega_3
		- I_{32} \omega_2 \omega_2
		- I_{31} \omega_1 \omega_2
		+ I_{21} \omega_1 \omega_3
		\\
		(I_{11} - I_{33})\omega_3 \omega_1 
		+ I_{31} \omega_1 \omega_1
		- I_{13} \omega_3 \omega_3
		- I_{12} \omega_2 \omega_3
		+ I_{32} \omega_2 \omega_1
		\\
		(I_{22} - I_{11})\omega_1 \omega_2 
		+ I_{12} \omega_2 \omega_2
		- I_{21} \omega_1 \omega_1
		- I_{23} \omega_3 \omega_1
		+ I_{13} \omega_3 \omega_2
		\\
		\end{pmatrix} \\
		I &= \begin{pmatrix}
		 I_{11} & -I_{12} & -I_{13}\\
		-I_{21} &  I_{22} & -I_{23}\\
		-I_{31} & -I_{32} &  I_{33}\\
		\end{pmatrix}
	\end{split}
\end{equation}

We use the equation (\ref{eq29}) to compute the values of $\omega_i(t+\Delta t/2)$ and $\dot \omega_i(t+\Delta t)$. Following the approach suggested in \cite{pfc2008} we use the iterative algorithm to find these unknowns:

\begin{itemize}
    \item Set $n=0$ \\
    \item Set $\omega_i^{[0]}$ to the initial angular velocity.
    \item (*) Solve (\ref{eq29}) for $\dot \omega_i$
    \item Determine a new (intermediate) angular velocity: $\omega_i^{[new]} = \omega_i^{[0]} + \dot \omega_i^{[n]} \Delta t$
    \item Revise the estimate of $\omega_i$ as:
    $\omega_i^{[n+1]} = 0.5 (\omega_i^{[0]} + \omega_i^{[new]})$
    \item Set $n:=n+1$ and go to (*)
\end{itemize}

This algorithm gives us the value of the angular velocity that is further used to update the position at the second step of leap-frog algorithm. The number of steps necessary for the sufficient precision varies depending on the application and is usually chosen in range of $2-5$.

The described approach is rather general, which potentially allows extension of the notion of clumps on quite wide set of pebble entities, including particles that do not track their orientations \cite{pfc2008}. However, the algorithm is inferior in terms of precision and performance compared to modern rigid-body integrators \cite{johnson2008, omelyan1998}, because of significant overhead related to solving the equations of motion in inertial frame -- this can be significant for clumps consisting of small numbers of pebbles, when the duration of rigid body integration is non-negligible compared to duration of updating positions of pebbles.

\subsection{Interaction of clump particles with periodic boundaries}

The complete description of the logic of interaction of spherical particles (classes \texttt{BaseParticle},
\texttt{SphericalParticle}) and periodic boundaries can be found in \cite{Lantman2019}. This logic had to be adjusted for rigid clumps. Below we briefly describe the corresponding modifications.

The original scheme utilizes the concept of primary particles and ``ghost'' particles that are introduced to represent interactions across periodic boundaries. ``Ghost'' particles are created when the primary particle approaches closely the periodic boundary, and ``switch'' status with the primary particle when the migration over the periodic boundary occurs. Our implementation introduced two minor modifications to this scheme to ensure correct treatment of rigid clumps in a periodic box:

\begin{itemize}
    \item ``Clump'' particles are never erased/created in a course of the simulation. They migrate over the periodic boundary seamlessly by direct specification of the position property. This way the necessity of sending the ``pebble'' pointers between ``clumps'' is avoided.

    \item The ``ghost'' particles for ``clumps'' do not exist, since no interaction is treated at the level of ``clump'' particles.

    \item The procedures of adding moments to a ``clump'' particle from the forces/moments acting on ``pebbles'', and computation of the translational velocity/position of ``pebble'' particles, utilize the minimum image convention to determine the length of the ``lever'' -- the vector connecting the center of ``clump'' particle (clump's COM), and the center of ``pebble'' particle. 
\end{itemize}

These adjustments are introduced in \texttt{/Drivers/Clump/ClumpHeaders/Mercury3DClump.h}, provide full functionality of all types of periodic boundaries, implemented earlier in \textit{MercuryDPM}.

\subsection{Random generation of non-overlapping clumps}

It is often necessary to create rigid clumps with random initial orientation. In order to provide equal probability of every orientation, we use the following scheme of clump random rotation: we first rotate the clump instance counterclockwise about $n_3$ direction by the angle $\alpha$, and then rotate the clump to match its principal direction $n_3$ with the random vector on a unit sphere $(\theta, \phi)$ in a spherical coordinate system: $n_3^{rot} = (\sin{\theta} \cos{\phi}, \sin{\theta} \sin{\phi}, \cos{\theta})$. The random values of $\alpha, \phi$ are chosen uniformly in the range $(0, 2 \pi)$, while the angle $\theta$ is chosen as $\arccos(p)$, where $p$ is uniformly distributed in $(-1,1)$. Such choice of random orientation angles ensures equal probability of every possible clump orientation.   

In order to ensure a placement of a new clump into the deposition domain without overlaps with the previous clumps, a straightforward algorithm is used to ensure that neither pebble of newly deposited clump overlaps with any pebble of the existing clump.

\subsection{Modifications of energy computing routines}

The routines computing rotational and translational kinetic energy of the clump, as well as its potential gravitational energy, had to be straightforwardly adjusted to reflect the correct inertial/gravitational properties of a clump, computed as detailed above.


\section{Examples}

\subsection{Computation of TOI -- precision of the summation}

This brief example illustrates the precision of our approaches used to compute mass and tensor of inertia of the clumps. The test model consists of two spherical pebbles of unit radius, with centers separated by one diameter of a pebble (Fig. 4 (A) gives the model represented with pebbles, tetrahedrons and voxels). This simple model allows immediate exact evaluation of inertial properties of this non-spherical, non-convex shape. The mass of the clump, as well as its major and minor moments of inertia are then evaluated with tetrahedral and voxel discretization. The vertices of tetrahedrons are the origin $(0,0,0)$ and the triangles constructed by equispaced angular subdivision of each pebble sphere on $N$ equal segments along latitude angle $\theta \in (0, \pi)$ and on $2N$ segments along azimuth angle $\phi \in (0, 2\pi)$ (see Fig.4 (A)). For the voxels, the refinement degree $N$ is defined as the number of voxels along the diameter of a pebble.

\begin{figure}
	\begin{center}
		\includegraphics[width=12 cm]{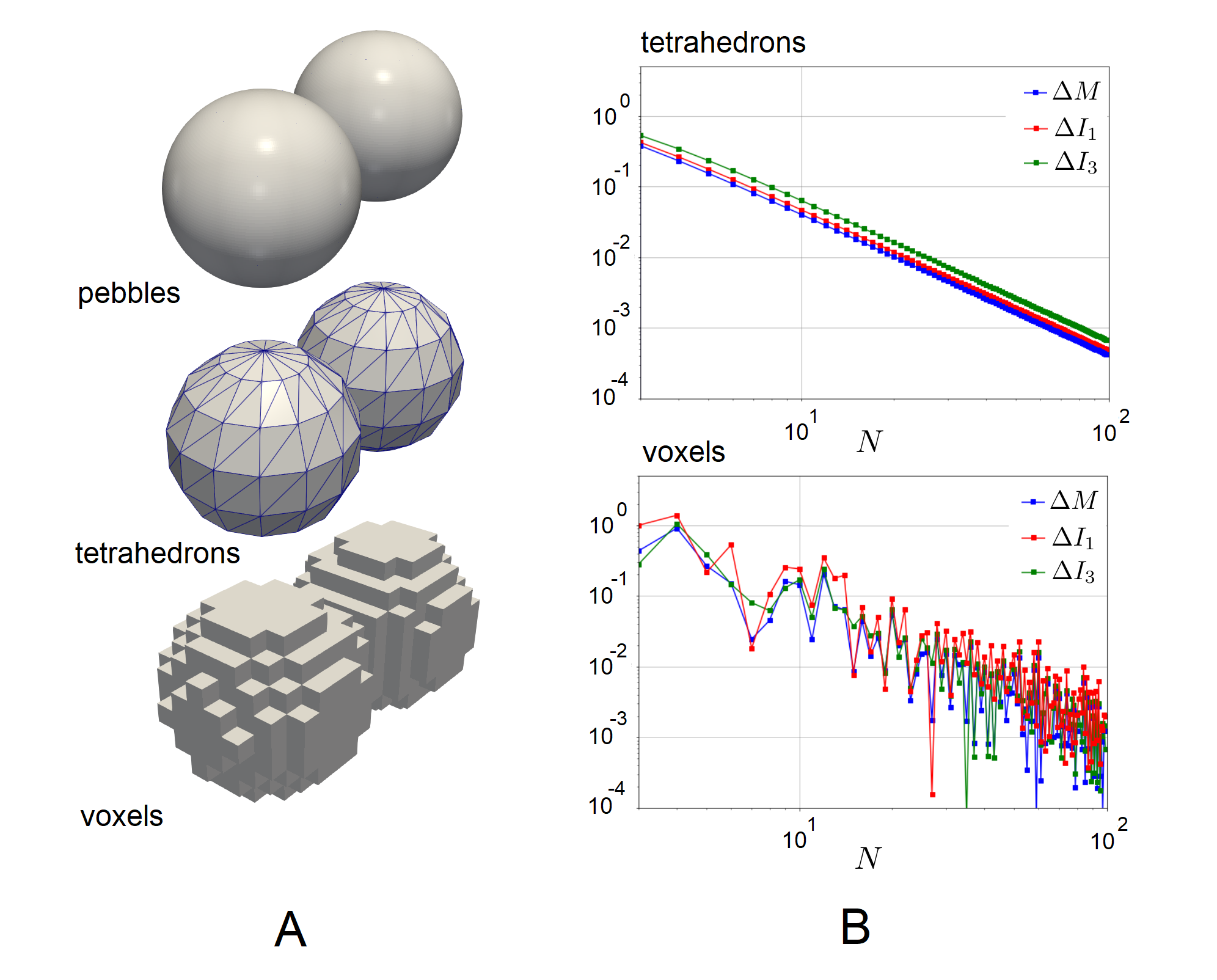}
		\protect\caption{(A) The model of a clump under test, represented with spherical pebbles, tetrahedrons and voxels. (B) Relative error in computation of clump's mass ($\Delta M$), major ($\Delta I_1$) and minor ($\Delta I_3$) principal components of inertia, as a function of the model refinement $N$ (see the definitions above), for tetrahedrons (top) and voxels (bottom).}
	\end{center}
\end{figure}

Fig. 4(B) demonstrates the convergence of relative error in computation of mass and principal moments of inertia with the degree of model refinement $N$. We can clearly see that the error is inversely proportional to $N$, both for tetrahedron and voxel discretization. The latter, however, features significant chaotic error, which suggests necessity of further improvement of an algorithm.

\subsection{Dynamics of a single particle - energy equipartition}

The simple simulation depicted in Fig. 5(A) is located at \texttt{Drivers/Clump/Single/Single.cpp}. An elastic, rod-like particle is placed into a cubic box  with elastic walls (no friction, no dissipation, linear contact model is employed). At the initial moment of simulation, the particle is assigned the initial translational velocity $V$, orientation along x axis and zero initial angular velocity $w$. After few collisions, the alignment of the particle with x axis breaks, and each next collision causes redistribution of energy between translational and rotational degrees of freedom (Fig. 5(B)). In a long enough timeline we see the energy equipartition between available degrees of freedom. For example, if the particle bounces strictly along $y$ axis between two elastic walls, and rotates around its principal axis co-oriented with $z$, it has only one translational and one rotational degree of freedom. We can therefore foresee that the equipartition will manifest itself with the ratio of 1 between the tranlational kinetic energy $mv^2/2$ and rotational kinetic energy $I\omega^2/2$ in a sufficiently long simulation. This is precisely what happens (Fig. 5(C)). Similarly, the different initial conditions leading to a different system of available degrees of freedom lead to different ratios. For example, if the initial translational velocity has two components, leading to two translational degrees of freedom, the ratio of rotational and translational energy converges to 0.5.

\subsection{Dynamics of a single particle - Dzhanibekov effect}

The example \texttt{Drivers/Clump/TBar/TBar.cpp} demonstrates so-called Dzhanibekov effect - instability of rotation around the second principal axis (see, e.g., \cite{Ashbaugh1991}). It manifests itself in a series of flips of an object rotating around its intermediate axis -- the classical example is a wingnut rotating around its axis in the condition of zero gravity. The simulation in this example reproduces this effect for  T-shaped clump (Fig. 6(A)), rotating around its second principal axis (see Video 1 in the supplementary information \cite{SimulatedDzhanibekov}). It is important to note that the observed angular momentum and rotational kinetic energy are well preserved during the simulation -- for example, as can be seen in Fig. 6(B), the relative drift of the rotational energy does not exceed $10^{-3}$ for $8$ flip cycles. 

\begin{figure}
	\begin{center}
		\includegraphics[width=16 cm]{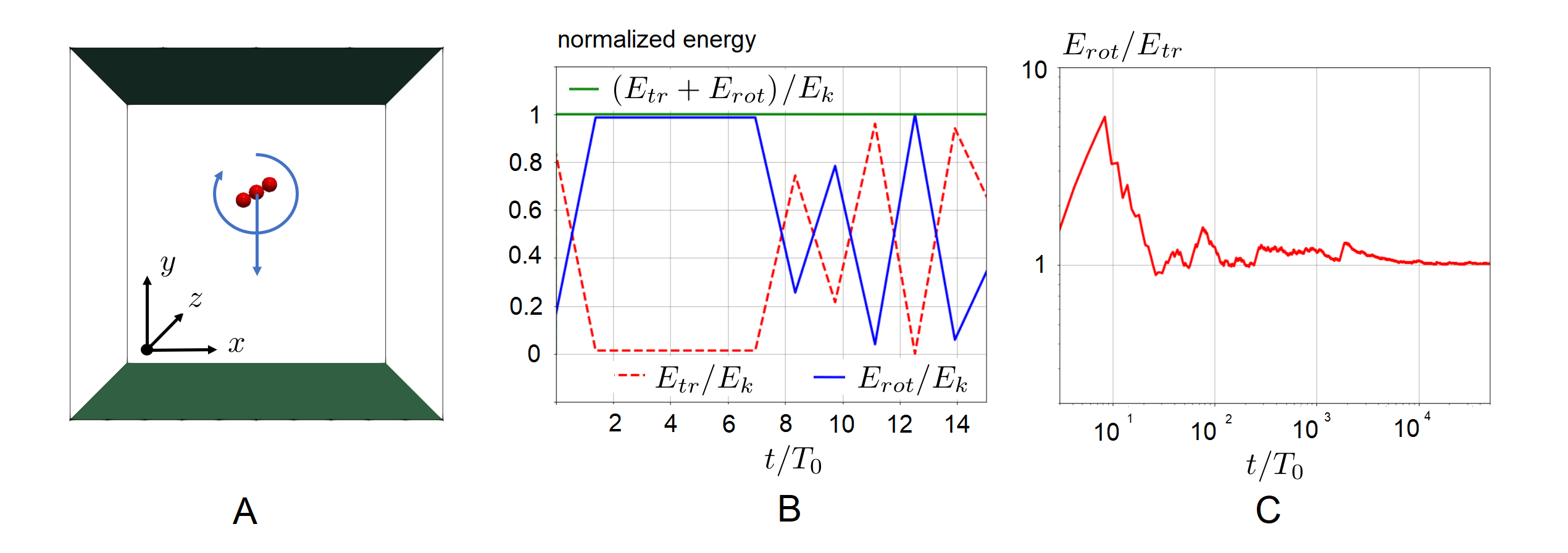}
		\protect\caption{(A) The model of a single-atom ideal gas with one translational and one rotational degree of freedom. (B) Observed fractions of translational and rotational kinetic energies as functions of time, for a time span comprising first 20 collisions. (C) The ratio between the rotational and translational kinetic energy, averaged over sufficiently long simulation time ($~5 \times 10^4$ particle-wall collisions).}
	\end{center}
\end{figure}

\begin{figure}
	\begin{center}
		\includegraphics[width=16 cm]{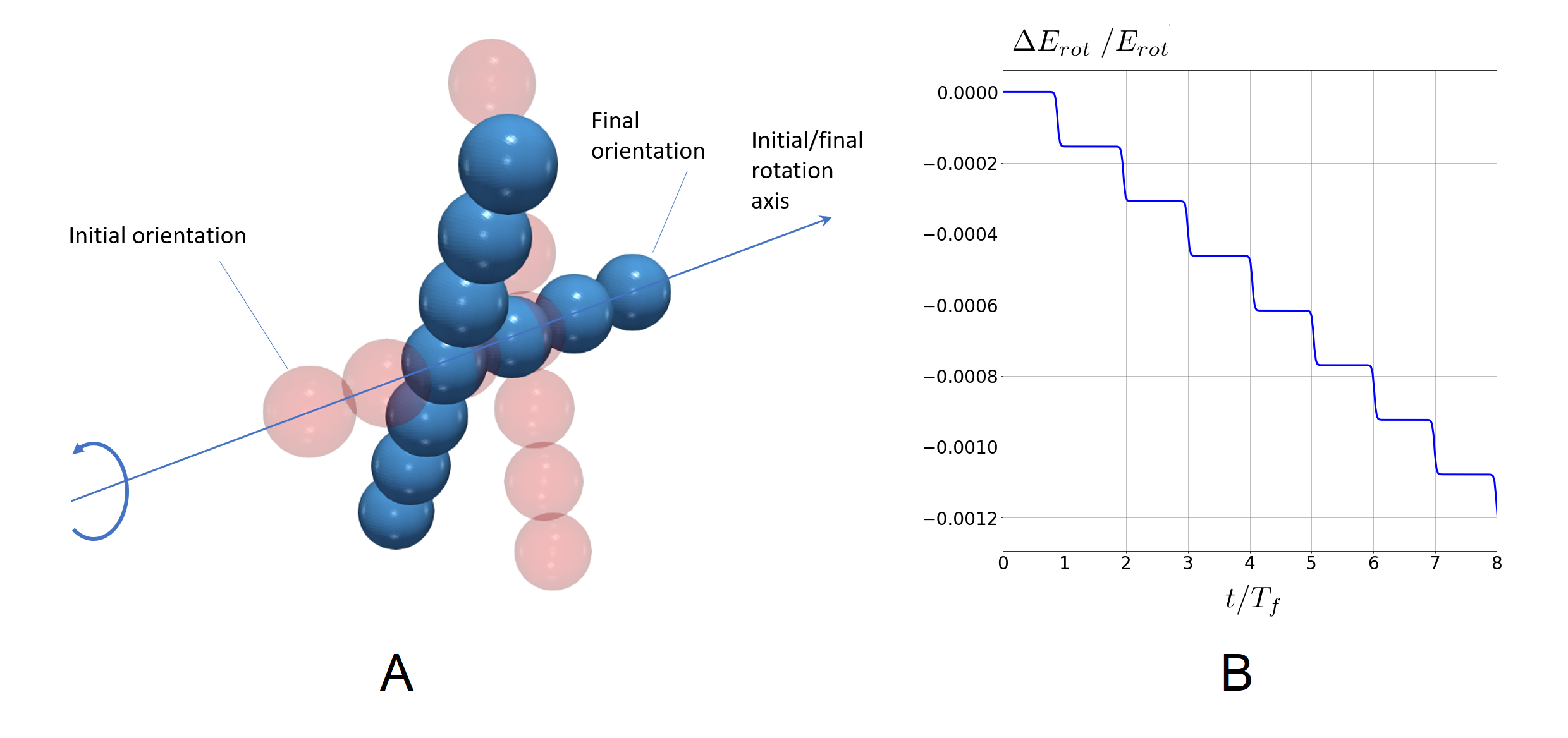}
		\protect\caption{(A) Evolution of the orientation of a T-bar, (B) observed relative drift of its kinetic energy.}
   
	\end{center}
\end{figure}

\subsection{Rolling of a Gömböc}

Gömböc is the convex body that, being put on the flat surface, has one point of stable and one point of unstable equilibrium \cite{Gomboc}. Arbitrarily oriented at the initial moment, provided sufficient energy dissipation, the gömböc finally arrives to its only stable equilibrium position. We use the model of a gömböc depicted in Fig. 7(A) to create a clump (Fig. 7(B)), mimicking the behavior of a gömböc. The clump was generated using the algorithm \cite{ferellec2010method} and has $182$ pebbles. We simulated the dynamics of gömböc shape, dropped to the flat surface (\texttt{./Drivers/Clump/Gomboc/Gomboc.cpp}) Our simulation (Video 2 in the supplementary imformation \cite{GombocVideo}) indicate that, after a series of metastable rotational oscillations (Fig. 7(C)), gömböc shape does indeed arrive to a unique stable orientation. Our experiments indicate that if the initial energy of a gömböc is too low, it may get stuck in one of the local energy minima that emerge due to approximation of the original shape by a finite number of spherical particles. Besides this effect, our simulations compare nicely with the experiments with real Gömböc shape.

\begin{figure}
	\begin{center}
		\includegraphics[width=16cm]{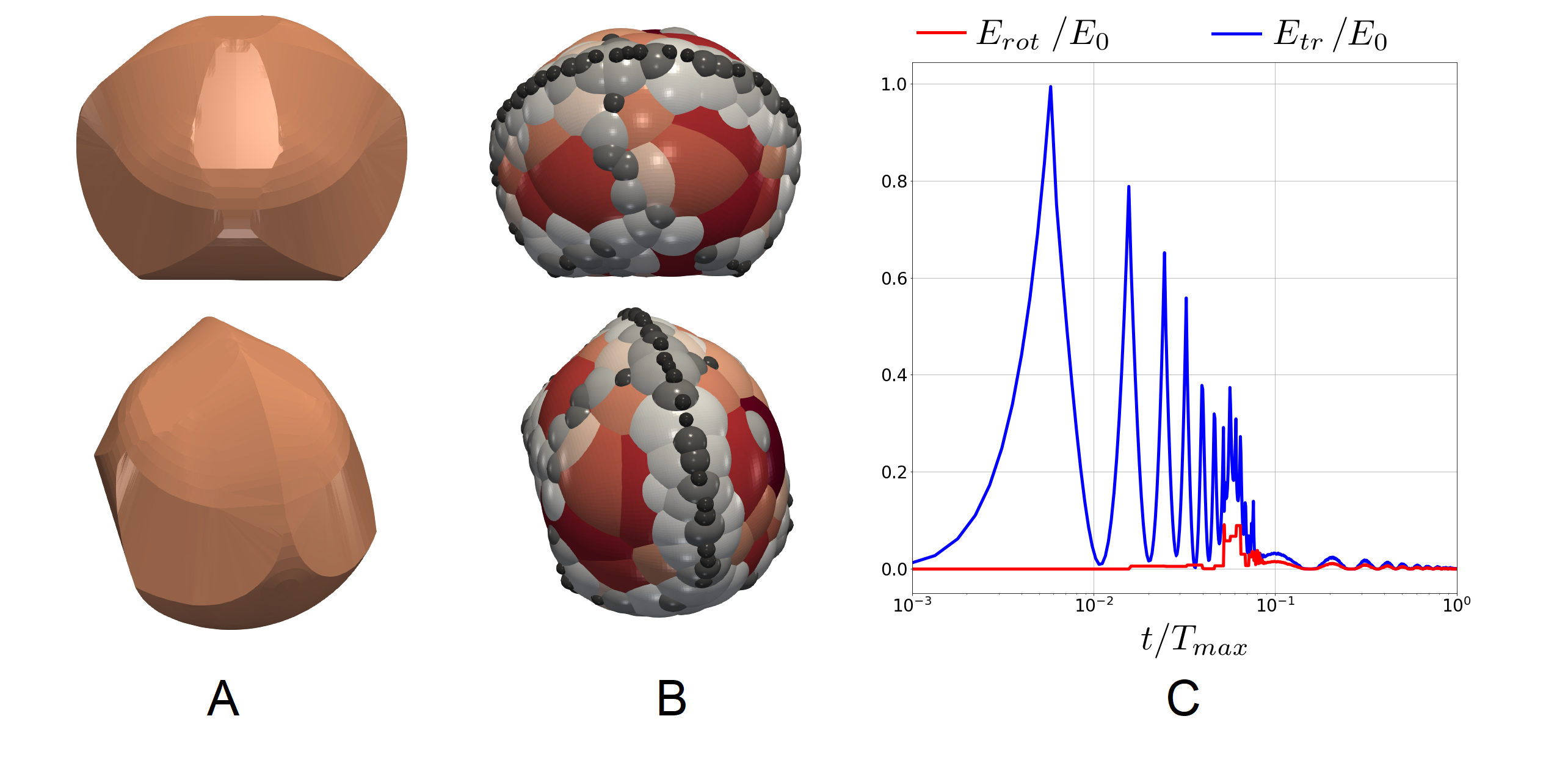}
		\protect\caption{Gömböc - (A) original stl model, (B) its rigid clump representation, computed according to \cite{ferellec2010method}, (C) Evolution of the translational and rotational kinetic energy with time in the simulation. The simulation duration was chosen to feature entire motion trajectory of a gömböc with realistic damping parameters.}
	\end{center}
\end{figure}

\subsection{Domino effect}

Domino effect is well know to be quite non-trivial benchmark example for DEM simulation with nonspherical particles \cite{Domino}. We provide a driver file designed for parametric studies of a domino effect (see \texttt{./Drivers/Clump/Domino/Domino.cpp}). Dominoes are rectangular regular packings of pebbles, equispaced along the straight line (Fig. 8(A)). At the initial moment the domino 1 is given an initial push with the cue - a spherical particle. The initial propagation of the domino wave is to a large extent affected by the initial velocity of the cue, however, the steady state velocity does not depend at all on this initial velocity. This, in particular, manifests itself in a constant time derivative of the potential energy (Fig. 8(B)) that does not depend on the initial cue velocity. This invariance of the domino wave velocity is well-known and often attributed \cite{Leeuwen_Domino} to dissipative effects; however, there are theoretical/numerical evidence \cite{EJ_Domino_2007, Domino} that it takes place even in the case of perfectly elastic collisions between the dominoes.

\begin{figure}
	\begin{center}
		\includegraphics[width=14.5cm]{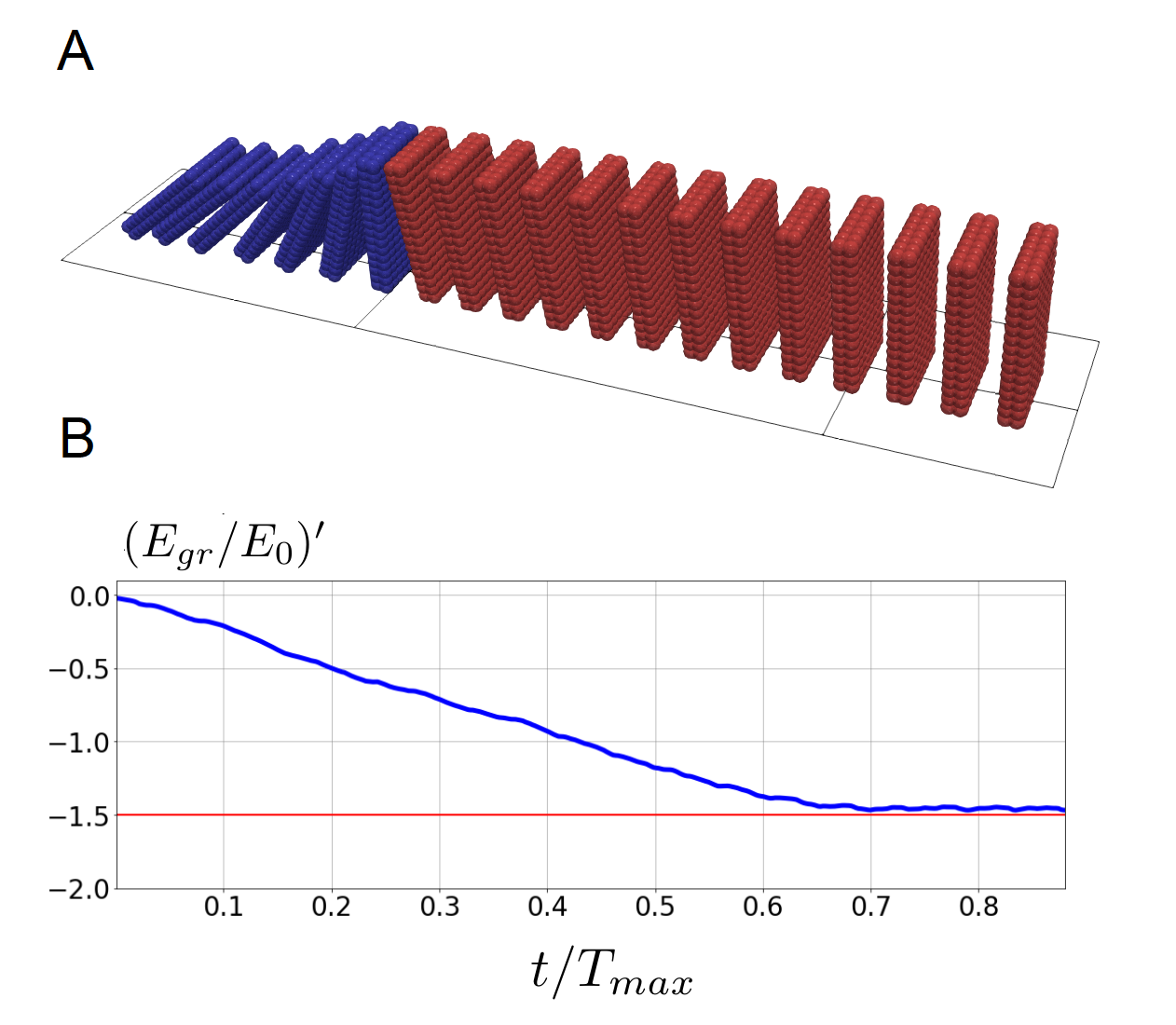}
		\protect\caption{(A) Geometry of DEM model of domino wave propagation, (B) Constant rate of change of the potential energy with time in the steady-state domino wave propagation ($E_0$ is the initial gravitational potential energy of the system; the simulation duration roughly corresponds to duration of domino wave propagation over 20 dominoes).}
	\end{center}
\end{figure}

\subsection{Dense gas of interacting T-shaped particles in a periodic box}

The driver $\texttt{Drivers/Clump/TGas/TGas.cpp}$ demonstrates the evolution of six hundreds of T-shaped rigid particles of arbitrary initial velocities, angular velocities and orientations, that are deposited in a triple periodic box without initial overlaps, with zero initial rotational velocities and random initial translational velocities (Fig. 9(A)). Shortly after the beginning of the simulation, we can see the complete energy equipartition (Fig. 9(B)). The driver code can be easily adjusted to introduce elastic walls, gravity, dissipation etc.  

\begin{figure}
	\begin{center}
		\includegraphics[width=16.0cm]{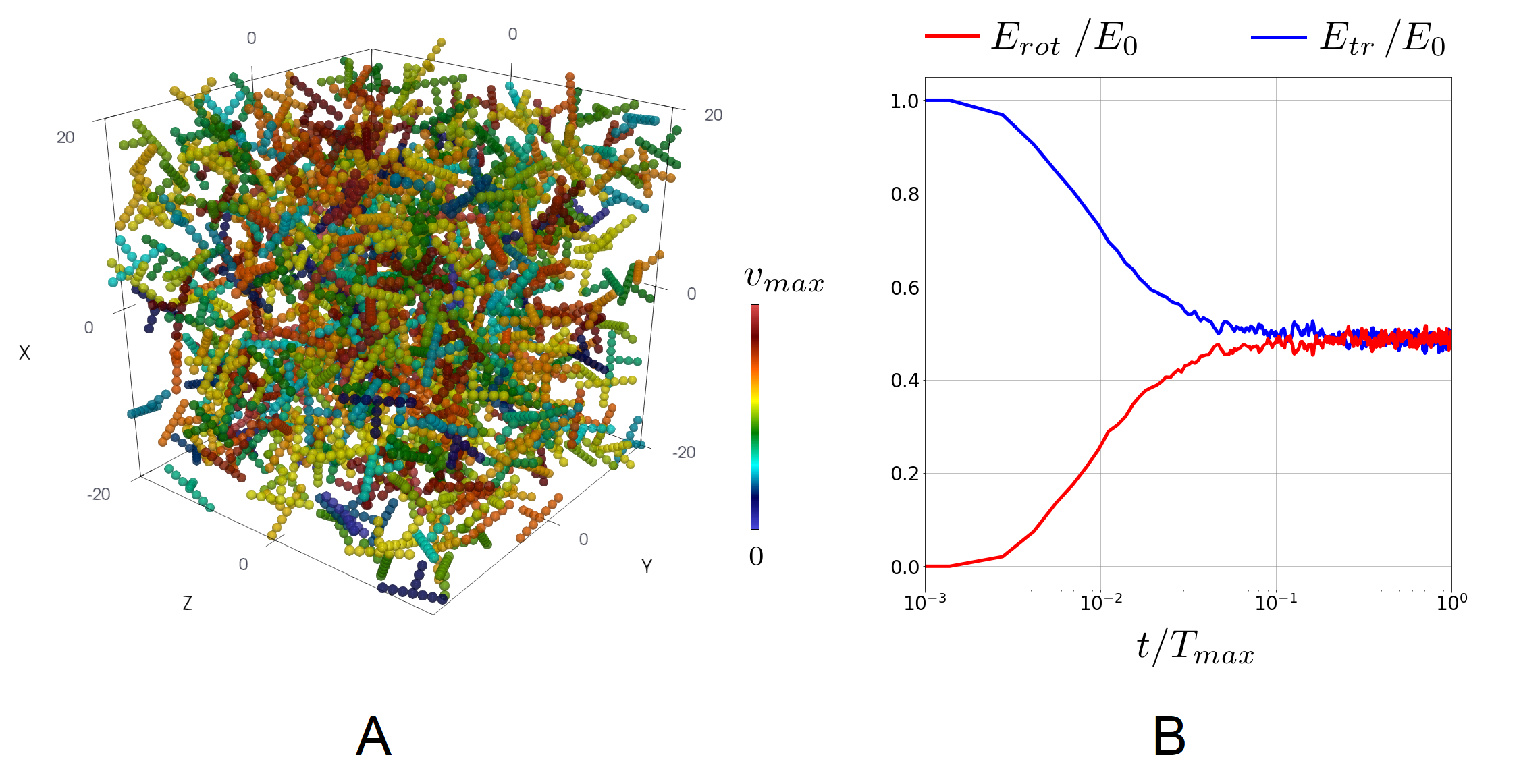}
		\protect\caption{Multiple T-bars in a box.(A) Initial geometry, (B) evolution of rotational and translational kinetic energy with time; the simulation duration was chosen to resolve the energy equipartition process.}
	\end{center}
\end{figure}

\subsection{Multiple clumps in a rotating drum}

A concluding example $\texttt{Drivers/Clump/RotatingDrum/RotatingDrum.cpp}$ features a collective motion of complex-shaped clumps in a rotating horizontal drum in the field of gravity (Fig. 10(A)). The gömböc shape described above was used as a clump instance, $27$ clumps were deposited in a volume of a drum without initial overlaps between themselves and the walls of the drum. The contact friction at both wall-clump and clump-clump contacts has zero rolling friction and high sliding friction of $0.6$. At the initial moment of simulation the drum starts to rotate with the constant angular velocity. The Video 3 in the supplementary information \cite{RotatingDrum} highlights the dynamic evolution of the system. Fig. 10(A) shows the geometry of the system, Fig. 10(B) gives the evolution of the gravitational potential energy of the clumps (normalized on the lowest energy observed in the beginning of the simulation) with time. One can see discrete events of sliding/repose of the bed ($8$ per $2$ full revolutions of the drum). This simulation validates the efficiency of the clump implementation in a moderate-size single-core simulation.

\begin{figure}
	\begin{center}
		\includegraphics[width=16.0cm]{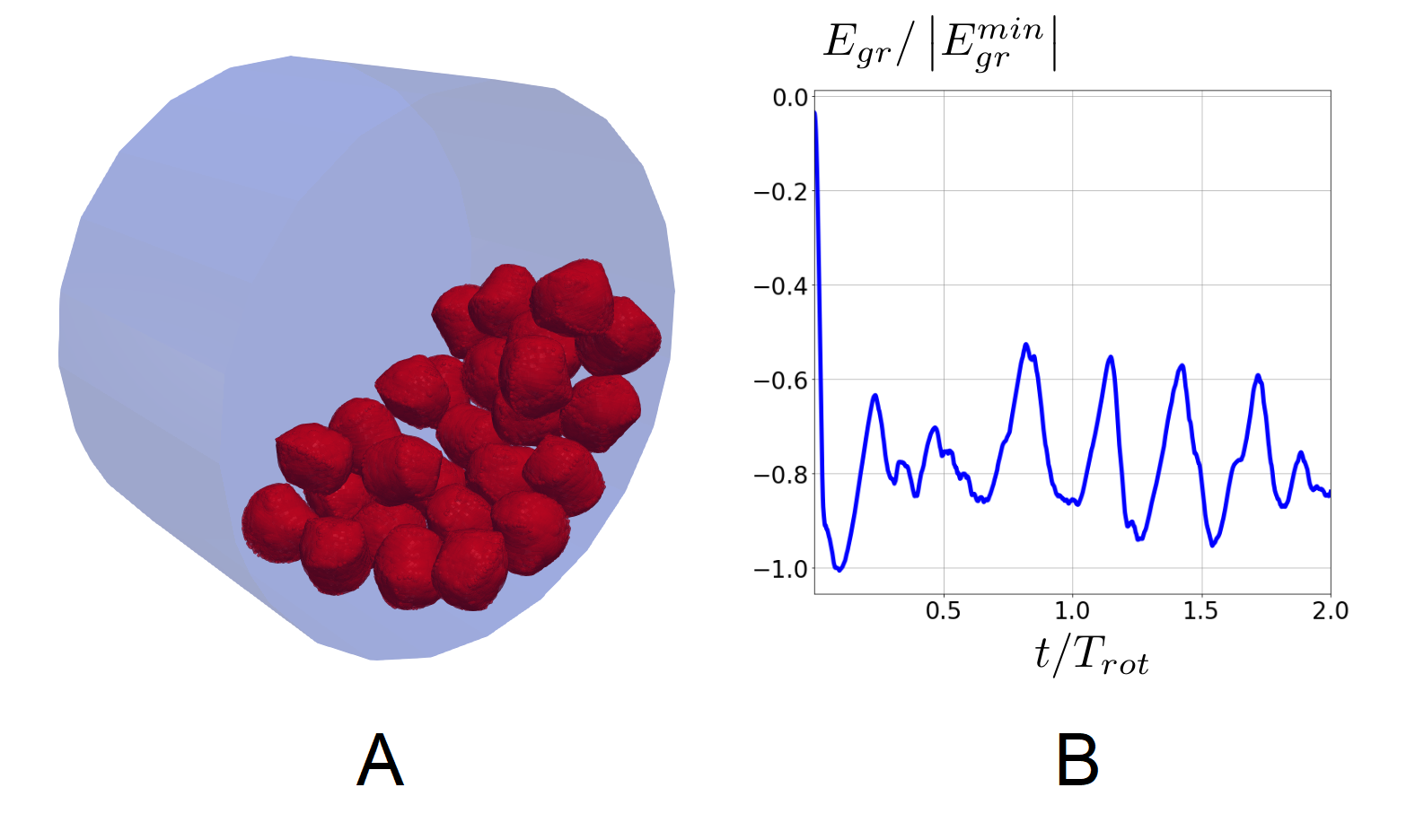}
		\protect\caption{Clumps in a rotating drum. (A) Problem geometry, (B) normalized potential energy of the clumps versus time, featuring sloshing motion pattern.}
	\end{center}
\end{figure}

\subsection{Efficiency of hierarchical grid contact detection algorithm for highly polydisperse clump systems}

One of the strong features of MercuryDPM is its efficient contact detection algorithm oriented on highly polydisperse particle assemblies \cite{Ogarko2012}. It is interesting to see how the single-core simulation performance of polydisperse clump systems is affected by the maximum number of levels of hierarchical grid employed by the contact detection algorithm (see \cite{Ogarko2012, Mercury2020} for details). Our benchmark examples predictably demonstrate that small models do not benefit from multiple levels of hierarchical grid used in contact detection, while larger models perform much faster with hierarchical grid turned on. The rotating drum simulation described above is used here to demonstrate the effect of multiple levels of hierarchical grid on the performance of simulation of the polydisperse clumps. Two (otherwise identical) simulations with different clump resolution were studied: \textbf{model 1} had clumps of $182$ pebbles ($4 914$ pebbles in total) and the size ratio of the largest to the smallest pebble of $28.83$; \textbf{model 2} had the same clump surface represented by $423$ pebbles ($11 421$ pebbles in total) with the size ratio of the largest to the smallest pebble of $53.36$ (Fig. 11(A)).  

\begin{figure}
	\begin{center}
		\includegraphics[width=16.0cm]{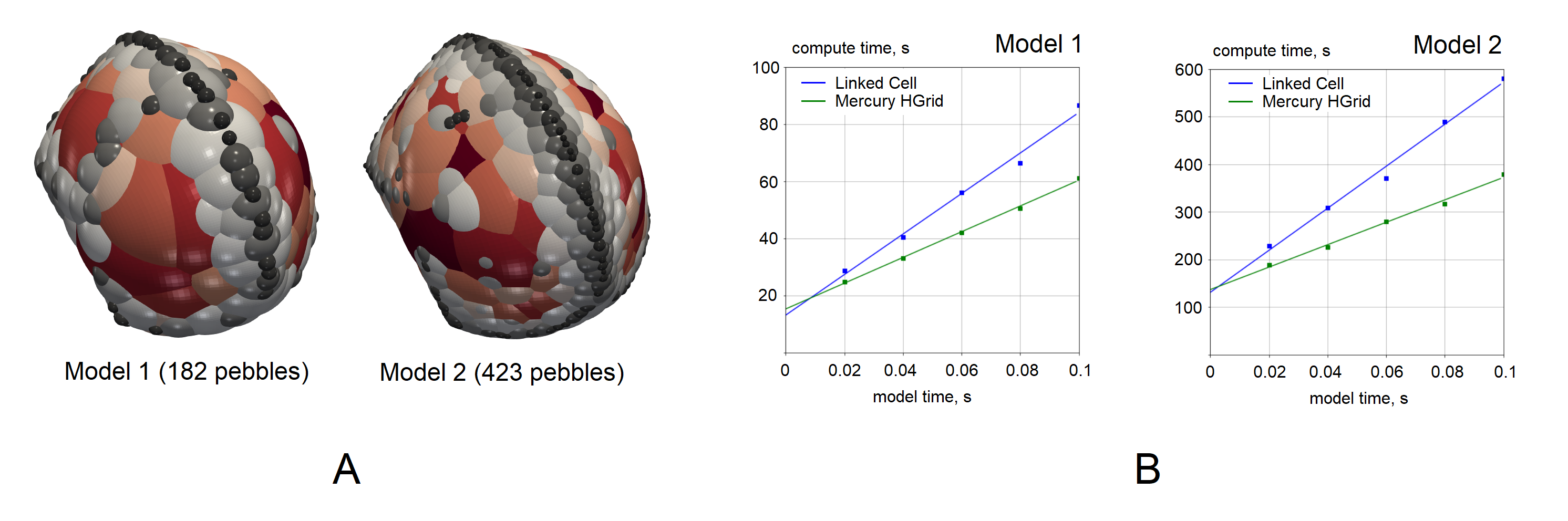}
		\protect\caption{Illustration of effect of multiple grid levels on computational performance.(A) Clumps used in model 1 and model2, (B) Dependence of compute time and model time for two models and contact detection approach used.}
	\end{center}
\end{figure}

Both models were studied in simulations with contact detection algorithm limited to one hierarchical grid level (regular linked cell algorithm, blue plots on Fig. 11(B)) and with three hierarchical grid levels (MercuryDPM default value, green plots on Fig. 11(B)). Accurate comparison of performance results in $57\%$ increase in the cycle-time performance for the \textbf{model 1} and $87\%$ - for the \textbf{model 2}. For larger models this increase in performance is expected to be even more dramatic \cite{Ogarko2012}. Therefore, we can see that MercuryDPM contact detection algorithm makes it well-suited for modeling polydisperse clumped particle systems.

\section{Conclusions}

This work details the implementation of rigid clumps within \textit{MercuryDPM} particle dynamics code. Necessary pre-processing tools, kernel modifications and driver files illustrating the applications are described. Due to advanced contact detection algorithm of \textit{MercuryDPM}, our implementation demonstrates high single-core performance for highly polydisperse clumps. The new features will certainly be useful to the \textit{MercuryDPM} community. The codes are currently available in the Master branch of the \textit{MercuryDPM} project \cite{MercuryMaster}. The implementation is under ongoing development, the changes in the existing implementation will be highlighted in the future release notes and the corresponding papers.  

\section*{Funding acknowledgements} 

\textit{MercuryDPM} has been supported by many projects, both past and present. The features presented here  were (partially) funded by the Dutch Research Council (NWO), in the framework of the ENW PPP Fund for the topsectors and from the Ministry of Economic Affairs in the framework of the “PPS-Toeslagregeling”.
 
\bibliographystyle{unsrtnat}
\bibliography{manuscript}

\end{document}